\documentclass[12pt]{article}
\usepackage[a4paper, margin=2.3cm]{geometry}
\usepackage{amsmath,amssymb,amsthm}
\usepackage{color}
\usepackage{parskip}

\newtheorem{theorem}{Theorem}[section]

\newtheorem{corollary}[theorem]{Corollary}

\newtheorem*{definition*}{Definition}

\def\F{\mathbb{F}}

\def\RR{\mathcal{R}}
\def\S{\mathcal{S}}
\newcommand{\I}{\mathcal{I}}
\begin{document}
\title{Conditional expanding bounds for two-variable functions over arbitrary fields}

\author{Hossein Nassajian Mojarrad\and 
    Thang Pham}
    \date{}
\maketitle
\begin{abstract}
In this short note, we use Rudnev's point-plane incidence bound to improve some results on conditional expanding bounds for two-variable functions over arbitrary fields due to Hegyv\'{a}ri and Hennecart \cite{heg}. 
\end{abstract}
\section{Introduction}
Throughout this chapter, by $\F$ we refer to any arbitrary field, while by $\mathbb{F}_p$, we only refer to the fields of prime order $p$. We denote the set of non-zero elements by $\mathbb{F}^*$ and $\mathbb{F}^*_p$, respectively. Furthermore, we use the following convention: if the characteristic of $\F$ is positive, then we denote its characteristic by $p$; if the characteristic of $\F$ is zero, then we set $p=\infty$. So a term like $N <p^{5/8}$ is restrictive in positive characteristic, but vacuous for zero one.

For $A\subset \mathbb{F}$,  the sum and the product sets are defined as follows:
\[ A+A=\{a+a' : a, a' \in A\}, ~A\cdot A=\{a\cdot a' : a, a' \in A\}.\]
For $A\subset \mathbb{F}_p$, Bourgain, Katz and Tao (\cite{bkt}) proved that if $ p^\delta < |A| < p^{1-\delta}$ for some $\delta > 0$,  then we have
\[ \max \left\lbrace |A+A|, |A\cdot A|\right\rbrace \gg |A|^{1+\epsilon},\]
for some $\epsilon = \epsilon(\delta) > 0$. Here, and throughout, by $X \ll Y$ we mean that there exists the constant $C>0$ such that $X\le  CY$.

In a breakthrough paper \cite{RRS}, Roche-Newton, Rudnev, and Shkredov improved and generalized this result to arbitrary fields. More precisely, they showed that for $A\subset \mathbb{F}$,  the sum set and the product set satisfy
\[\max\left\lbrace |A\pm A|, |A\cdot A|\right\rbrace \gg |A|^{6/5}, ~\max\left\lbrace |A\pm A|, |A : A|\right\rbrace \gg |A|^{6/5}.\]
We note that the same bound also holds for $|A(1+A)|$ \cite{SZ}, and $|A+A^2|$, $\max\left\lbrace |A+A|, |A^2+A^2|\right\rbrace$ \cite{P1}. We refer the reader to \cite{AMRS, BT, RRS, M} and references therein for recent results on the sum-product topic. 

Let $G$ be a subgroup of $\mathbb{F}^*$, and $g\colon G\to \mathbb{F}^*$ be an arbitrary function. We define
\[\mu(g)=\max_{t\in \mathbb{F}^*}\left\vert \left\lbrace x\in G\colon g(x)=t\right\rbrace \right\vert.\]

For $A, B\subset \mathbb{F}_p$ and two-variable functions $f(x,y)$ and $g(x,y)$ in $\mathbb{F}_p[x, y]$,  Hegyv\'{a}ri and Hennecart \cite{heg}, using graph theoretic techniques,  proved that if $|A|=|B|=p^\alpha$, then \[\max\left\lbrace |f(A,B)|, |g(A,B)|\right\rbrace\gg |A|^{1+\Delta(\alpha)},\] for some $\Delta(\alpha)>0$. More precisely, they established the following results. 
\begin{theorem}[\textbf{Hegyv\'{a}ri and Hennecart}, \cite{heg}]\label{thm1}
Let $G$ be a subgroup of $\mathbb{F}_p^*$. Consider the function $f(x,y)=g(x)(h(x)+y)$ on $G\times \mathbb{F}_p^*$, where $g,h\colon G\to \mathbb{F}_p^*$ are arbitrary functions. Define $m=\mu(g\cdot h)$. For any subsets $A\subset G$ and $B,C\subset \mathbb{F}_p^*$, we have
\[\left\vert f(A,B)\right\vert \left\vert B\cdot C\right\vert\gg \min\left\lbrace\frac{|A||B|^2|C|}{pm^2}, \frac{p|B|}{m}\right\rbrace.\] 
\end{theorem}
\begin{theorem}[\textbf{Hegyv\'{a}ri and Hennecart}, \cite{heg}]\label{2.2}
Let $G$ be a subgroup of $\mathbb{F}_p^*$. Consider the function $f(x,y)=g(x)(h(x)+y)$ on $G\times \mathbb{F}_p^*$, where $g,h\colon G\to \mathbb{F}_p^*$ are arbitrary functions. Define $m=\mu(g)$. For any subsets $A\subset G$, $B,C \subset \mathbb{F}_p^*$, we have
\[|f(A,B)||B+C|\gg \min\left\lbrace  \frac{|A||B|^2|C|}{pm^2}, \frac{p|B|}{m}\right\rbrace\]
\end{theorem}
It is worth noting that Theorem $6$ established by Bukh and Tsimerman  \cite{BT} does not cover such a function defined in Theorem \ref{2.2}. The reader can also find the generalizations of Theorems \ref{thm1} and \ref{2.2} in the setting of finite valuation rings in \cite{Ham}.

Suppose $f(x,y)=g(x)(h(x)+y)$ with $\mu(g), \mu(h)=O(1)$ and $A=B=C$. Then, it follows from Theorems \ref{thm1} and \ref{2.2} that 
\begin{enumerate}
\item If $|A|\gg p^{2/3}$, then we have
\[|f(A, A)||A\cdot A|,|f(A, A)||A+A|\gg p|A|.\]
\item If $|A|\ll p^{2/3}$, then we have
\begin{equation}\label{xxx1}|f(A, A)||A\cdot A|,|f(A, A)||A+A|\gg |A|^4/p.\end{equation}
\end{enumerate}
The main goal of this paper is to improve and generalize Theorems \ref{thm1} and \ref{2.2} to arbitrary fields for small sets. Our first result is an improvement of Theorem \ref{thm1}.
\begin{theorem}\label{thm1*}
Let $f(x,y)=g(x)(h(x)+y)$ be a function defined on $\mathbb{F}^*\times \mathbb{F}^*$, where $g,h\colon \mathbb{F}^*\to \mathbb{F}^*$ are arbitrary functions. Define $m=\mu(g\cdot h)$. For any subsets $A, B,C\subset \mathbb{F}^*$ with $|A|, |B|, |C|\le p^{5/8}$, we have
\[\max\left\lbrace |f(A, B)|, |B\cdot C|\right\rbrace\gg \min\left\lbrace\frac{|A|^{\frac{1}{5}}|B|^{\frac{4}{5}}|C|^{\frac{1}{5}}}{m^{\frac{4}{5}}}, \frac{|B||C|^{\frac{1}{2}}}{m}, \frac{|B||A|^{\frac{1}{2}}}{m}, \frac{|B|^{\frac{2}{3}}|C|^{\frac{1}{3}}|A|^{\frac{1}{3}}}{m^{\frac{2}{3}}}\right\rbrace.\]
\end{theorem}
The following are consequences of Theorem \ref{thm1*}.
\begin{corollary}\label{co1}
Let $f(x,y)=g(x)(h(x)+y)$ be a function defined on $\mathbb{F}^*\times \mathbb{F}^*$, where $g,h\colon \mathbb{F}^*\to \mathbb{F}^*$ are arbitrary functions with $\mu(g\cdot h)=O(1)$. For any subset $A\subset \mathbb{F}^*$ with $|A|\le p^{5/8}$, we have
\[\max\left\lbrace |f(A, A)|, |A\cdot A|\right\rbrace\gg |A|^{\frac{6}{5}}.\]
\end{corollary}
\begin{corollary}\label{corx}
Consider the subsets $A\subset \mathbb{F}^{*}$, and $ B, C\subset \mathbb{F}$ with $|A|, |B|, |C|\le p^{5/8}$. 
\begin{enumerate}
\item By fixing $g(x)=1$ and $h(x)=x^{-1}$, we get 
\[\max\left\lbrace |A^{-1}+B|, |B\cdot C|\right\rbrace\gg \min\left\lbrace|A|^{\frac{1}{5}}|B|^{\frac{4}{5}}|C|^{\frac{1}{5}}, |B||C|^{\frac{1}{2}}, |B||A|^{\frac{1}{2}}, |B|^{\frac{2}{3}}|C|^{\frac{1}{3}}|A|^{\frac{1}{3}}\right\rbrace.\]
\item By fixing $g(x)=x$ and $h(x)=1$, we have 
\[\max\left\lbrace |A(B+1)|, |B\cdot C|\right\rbrace\gg \min\left\lbrace|A|^{\frac{1}{5}}|B|^{\frac{4}{5}}|C|^{\frac{1}{5}}, |B||C|^{\frac{1}{2}}, |B||A|^{\frac{1}{2}}, |B|^{\frac{2}{3}}|C|^{\frac{1}{3}}|A|^{\frac{1}{3}}\right\rbrace.\]
\end{enumerate}
\end{corollary}
It follows from Corollary \ref{corx}(2) that if $B=A$ and $C=A+1$ then we have $|A(A+1)|\gg |A|^{6/5}$, which recovers the result of Stevens and de Zeeuw \cite{SZ}.

Our next result is the additive version of Theorem \ref{thm1*}, which improves Theorem \ref{2.2}.
\begin{theorem}\label{thm2*}
Let $f(x,y)=g(x)(h(x)+y)$ be a function defined on $\mathbb{F}^*\times \mathbb{F}^*$, where $g\colon \mathbb{F}^*\to \mathbb{F}^*$ are arbitrary functions. Define $m=\mu(g)$. For any subsets $A, B,C\subset \mathbb{F}^*$ with $|A|, |B|, |C|\le p^{5/8}$, we have
\[\max\left\lbrace |f(A, B)|, |B+ C|\right\rbrace\gg \min\left\lbrace\frac{|A|^{\frac{1}{5}}|B|^{\frac{4}{5}}|C|^{\frac{1}{5}}}{m^{\frac{4}{5}}}, \frac{|B||C|^{\frac{1}{2}}}{m}, \frac{|B||A|^{\frac{1}{2}}}{m}, \frac{|B|^{\frac{2}{3}}|C|^{\frac{1}{3}}|A|^{\frac{1}{3}}}{m^{\frac{2}{3}}}\right\rbrace.\]
\end{theorem}
\begin{corollary}\label{co2}
Let $f(x,y)=g(x)(h(x)+y)$ be a function defined on $\mathbb{F}^*\times \mathbb{F}^*$, where $g\colon \mathbb{F}^*\to \mathbb{F}^*$ are arbitrary functions with $\mu(g)=O(1)$. For any subset $A\subset \mathbb{F}^*$ with $|A|\le p^{5/8}$, we have
\[\max\left\lbrace |f(A, A)|, |A+ A|\right\rbrace\gg |A|^{\frac{6}{5}}.\]
\end{corollary}
Let $g(x)=x$ and $h(x)=1$, we have the following corollary.
\begin{corollary}
For $A, B, C\subset \mathbb{F}$ with $|A|, |B|, |C|\le p^{5/8}$, we have 
\[\max\left\lbrace |A(B+1)|, |B+ C|\right\rbrace\gg \min\left\lbrace|A|^{\frac{1}{5}}|B|^{\frac{4}{5}}|C|^{\frac{1}{5}}, |B||C|^{\frac{1}{2}}, |B||A|^{\frac{1}{2}}, |B|^{\frac{2}{3}}|C|^{\frac{1}{3}}|A|^{\frac{1}{3}}\right\rbrace.\]
\end{corollary}
By fixing $g(x)=x$ and $h(x)=0$, we have the following result. 
\begin{corollary}
For $A, B, C\subset \mathbb{F}$ with $|A|, |B|, |C|\ll p^{5/8}$, we have 
\[\max\left\lbrace |A\cdot B|, |B+C|\right\rbrace \gg \min\left\lbrace|A|^{\frac{1}{5}}|B|^{\frac{4}{5}}|C|^{\frac{1}{5}}, |B||C|^{\frac{1}{2}}, |B||A|^{\frac{1}{2}}, |B|^{\frac{2}{3}}|C|^{\frac{1}{3}}|A|^{\frac{1}{3}}\right\rbrace\]
\end{corollary}
In the case $A=B=C$, we recover the following result due to Roche-Newton, Rudnev, and Shkredov \cite{RRS}, which says that $\max\left\lbrace|A+A|, |A\cdot A|\right\rbrace\gg |A|^{6/5}$.

It has been shown in \cite{SZ} that if $f(x,y)=x(x+y)$, then $|f(A, A)|\gg |A|^{5/4}$  under the condtion $|A|\le p^{2/3}$. In the following theorem, we show that if either $|A+A|$ or $|A\cdot A|$ is sufficiently small, the exponent $5/4$ can be improved from the polynomials to a larger family of functions on $\mathbb{F}^*\times \mathbb{F}^*$
\begin{theorem}\label{thmx}
Let $f(x,y)=g(x)(h(x)+y)$ be a function defined on $\mathbb{F}^*\times \mathbb{F}^*$, where $g,h\colon \mathbb{F}^*\to \mathbb{F}^*$ are arbitrary functions with $\mu(f), \mu(g)=O(1)$. Consider the subset $A\subset \mathbb{F}^*$ with $|A|\le p^{5/8}$, satisfying
\[\min\left\lbrace|A+A|, |A\cdot A|\right\rbrace\le |A|^{\frac{9}{8}-\epsilon}\]
for some $\epsilon>0$. Then, we have 
\[|f(A, A)|\gg |A|^{\frac{5}{4}+\frac{2\epsilon}{3}}.\]
\end{theorem}
\section{Proofs of Theorems \ref{thm1*}, \ref{thm2*}, and \ref{thmx}}
Let $\RR$ be a set of points in $\mathbb{F}^3$ and $\S$ be a set of planes in $\mathbb{F}^3$.  We write $\I(\RR,\S) = |\{(r,s)\in \RR\times \S : r\in s\}|$ for the number of \emph{incidences} between $\RR$ and $\S$. To prove Theorems \ref{thm1*} and \ref{thm2*}, we make use of the following point-plane incidence bound due to Rudnev \cite{R}. A short proof can be found in \cite{Z}.
\begin{theorem}[\textbf{Rudnev}, \cite{R}]\label{thm:rudnev}
Let $\RR$ be a set of points in $\mathbb{F}^3$ and let $\S$ be a set of planes in $\mathbb{F}^3$, with $|\RR|\ll |\S|$ and $|\RR|\ll p^2$.
Assume that there is no line containing $k$ points of $\RR$.
Then
\[ \I(\RR,\S)\ll |\RR|^{1/2}|\S| +k|\S|.\] 
\end{theorem} 
\paragraph{Proof of Theorem \ref{thm1*}:}
Define $f(A, B)=\{f(a, b)\colon a\in A, b\in B\}, g(A)=\{g(a)\colon a\in A\}$, $h(A)=\{h(a)\colon a\in A\}$.
For $\lambda\in B\cdot C$, let 
\[E_\lambda=\left\vert \left\lbrace \left(f(a, b), c\cdot g(a)^{-1}, c\cdot h(a)\right)\colon (a, b, c)\in A\times B\times C, ~f(a, b)\cdot c\cdot g(a)^{-1}-c\cdot h(a)=\lambda\right\rbrace\right\vert,\]
where by $g(a)^{-1}$ we mean the multiplicative inverse of $g(a)$ in $\mathbb{F}^*$. For a given triple $(x, y, z)\in (\mathbb{F}^*)^3$, we count the number of solutions $(a,b,c)\in A\times B \times C$ to the following system
\[g(a)(h(a)+b)=x, ~c\cdot g(a)^{-1}=y, ~c\cdot h(a)=z.\]
This implies that 
\[ g(a)h(a)=zy^{-1}.\]
Since $\mu(g\cdot h)=m$, there are at most $m$ different values of $a$ satisfying the equation $g(a)h(a)=zy^{-1}$, and $b, c$ are uniquely determined in term of $a$ by the first and second equations of the system. This implies that 
\[ |A||B||C|/m\le \sum_{\lambda\in B\cdot C}E_\lambda .\]
By the Cauchy-Schwarz inequality, we get 
\begin{equation}\label{cn1} \left(|A||B||C|/m\right)^2 \le \left(\sum_{\lambda\in B\cdot C}E_{\lambda}\right)^2\le E\cdot |B\cdot C| ,\end{equation}
where $E=\sum_{\lambda\in B\cdot C}E_{\lambda}^2$.

Define the point set $\mathcal{R}$ as
\[ \mathcal{R} = \left\lbrace\left(c\cdot g(a)^{-1}, c\cdot h(a), g(a')(h(a')+b')\right): a, a'\in A, b'\in B, c\in C\right\rbrace\]
and the set of planes $\mathcal{S}$ as
\[ \mathcal{S} = \left\lbrace g(a)(h(a)+b)X - Y -c'g(a')^{-1}Z = -c'\cdot h(a')\colon a, a'\in A, b\in B, c'\in C\right\rbrace.\]
We have $E \le I(\RR, \S)$, and $|\mathcal{R}|=|\mathcal{S}|\le |f(A, B)||A||C|$. To apply Theorem \ref{thm:rudnev}, we need to find an upper bound on $k$ which is the maximum number of collinear points in $\mathcal{R}$. The projection of $\RR$ into the first two coordinates is the set $\mathcal{T}=\{(c\cdot g(a)^{-1}, c\cdot h(a))\colon a\in A, c\in C\}$. The set $\mathcal{T}$ can be covered by the lines of the form $y=g(a)h(a)x$ with $a\in A$. This implies that $\mathcal{T}$ can be covered by at most $|A|$ lines passing through the origin, with each line containing $|C|$ points of $\mathcal{T}$. Therefore, a line in $\mathbb{F}^3$ contains at most $\max\{|A|,|C|\}$  points of $\RR$, unless it is vertical, in which case it contains at most $|f(A, B)|$ points. In other words, we get
\[k\le \max\{|A|, |C|, |f(A, B)|\}.\]
If $|\RR|\gg p^2$, then we get $|f(A, B)||A||C|\gg p^2$. Since $|A|, |C|\le p^{5/8}$, we have $|f(A, B)|\gg p^{3/4}\gg |A|^{\frac{1}{5}}|B|^{\frac{4}{5}}|C|^{\frac{1}{5}}$, and we are done in this case. Thus, we can assume that $|\RR|\ll p^2$. Applying Theorem \ref{thm:rudnev}, we obtain
\begin{equation}\label{cn2}I(\RR, \S)\le |f(A, B)|^{3/2}|A|^{3/2}|C|^{3/2}+k|f(A, B)|A||C|.\end{equation}
Putting (\ref{cn1}) and (\ref{cn2}) together gives us 
\[\max\left\lbrace |f(A, B)|, |B\cdot C|\right\rbrace\gg \min\left\lbrace\frac{|A|^{\frac{1}{5}}|B|^{\frac{4}{5}}|C|^{\frac{1}{5}}}{m^{\frac{4}{5}}}, \frac{|B||C|^{\frac{1}{2}}}{m}, \frac{|B||A|^{\frac{1}{2}}}{m},  \frac{|B|^{\frac{2}{3}}|C|^{\frac{1}{3}}|A|^{\frac{1}{3}}}{m^{\frac{2}{3}}}\right\rbrace.\]
This completes the proof of the theorem. $\square$
\paragraph{Proof of Theorem \ref{thm2*}:}
The proof goes in the same direction as Theorem \ref{thm1*}, but for the sake of completeness, we include the detailed proof. For $\lambda\in B+ C$, let 
\[E_\lambda=\left\vert \left\lbrace \left(f(a, b), g(a)^{-1}, c-h(a)\right)\colon (a, b, c)\in A\times B\times C, ~f(a, b)\cdot g(a)^{-1}+(c-h(a))=\lambda\right\rbrace\right\vert.\]
For a given triple $(x, y, z)\in (\mathbb{F}^*)^3$, we count the number of solutions $(a,b,c)\in A\times B \times C$ to the following system
\[g(a)(h(a)+b)=x, ~g(a)^{-1}=y, ~c-h(a)=z.\]
Since $\mu(g)=m$, there are at most $m$ different values of $a$ satisfying the equation $g(a)=y^{-1}$, and $b, c$ are uniquely determined in term of $a$ by the first and  third equations of the system.
This implies that 
\[|A||B||C|/m \le\sum_{\lambda\in B+C}E_\lambda .\]
By the Cauchy-Schwarz inequality, we have 
\begin{equation}\label{cn21} \left(|A||B||C|/m\right)^2 \le \left(\sum_{\lambda\in B+C}E_{\lambda}\right)^2\le E\cdot |B+C| ,\end{equation}
where $E=\sum_{\lambda\in B+C}E_{\lambda}^2$.
Define the point set $\mathcal{R}$ as
\[ \mathcal{R} = \left\lbrace\left(g(a)^{-1}, c-h(a), g(a')(h(a')+b')\right): a, a'\in A, b'\in B, c\in C\right\rbrace,\]
and the collection of planes $\mathcal{S}$ as
\[ \mathcal{S} = \left\lbrace g(a)(h(a)+b)X +Y -g(a')^{-1}Z = c'-h(a')\colon a, a'\in A, b\in B, c'\in C\right\rbrace.\]
It is clear that  $|\mathcal{R}|=|\mathcal{S}|\le |f(A, B)||A||C|$, and $E\le I(\mathcal{R}, \mathcal{S})$. To apply Theorem \ref{thm:rudnev}, we need to find an upper bound on $k$ which is the maximum number of collinear points in $\mathcal{R}$.  The projection of $\RR$ into the first two coordinates is the set $\mathcal{T}=\{\left(g(a)^{-1}, c-h(a)\right)\colon a\in A, c\in C\}$. The set $\mathcal{T}$ can be covered by at most $|A|$ lines of the form $x=g(a)^{-1}$ with $a\in A$, where each line contains $|C|$ points of $\mathcal{T}$. Therefore, a line in $\mathbb{F}^3$ contains at most $\max\{|A|,|C|\}$  points of $\RR$, unless it is vertical, in which case it contains at most $|f(A, B)|$ points. So we get
\[k\le \max\{|A|, |C|, |f(A, B)|\}.\]
If $|\RR|\gg p^2$, this implies that $|f(A, B)||A||C|\gg p^2$. Since $|A|, |C|\le p^{5/8}$, we have $|f(A, B)|\gg p^{3/4}\gg |A|^{\frac{1}{5}}|B|^{\frac{4}{5}}|C|^{\frac{1}{5}}$, and we are done. Thus, we can assume that $|\RR|\ll p^2$. Applying Theorem \ref{thm:rudnev}, we obtain
\begin{equation}\label{cn22}I(\RR, \S)\le |f(A, B)|^{3/2}|A|^{3/2}|C|^{3/2}+k|f(A, B)|A||C|.\end{equation}
Putting (\ref{cn21}) and (\ref{cn22}) together gives us 
\[\max\left\lbrace |f(A, B)|, |B+C|\right\rbrace\gg \min\left\lbrace\frac{|A|^{\frac{1}{5}}|B|^{\frac{4}{5}}|C|^{\frac{1}{5}}}{m^{\frac{4}{5}}}, \frac{|B||C|^{\frac{1}{2}}}{m}, \frac{|B||A|^{\frac{1}{2}}}{m}, \frac{|B|^{\frac{2}{3}}|C|^{\frac{1}{3}}|A|^{\frac{1}{3}}}{m^{\frac{2}{3}}}\right\rbrace.\]
This completes the proof. $\square$
\paragraph{Proof of Theorem \ref{thmx}:}
One can assume that $|f(A,A)| \le |A|^2$,  since otherwise we are done. Now by the proofs of Theorems \ref{thm1*} and \ref{thm2*} for $A\subset \mathbb{F}^*$ with $|A|\le p^{5/8}$, we have 
\[|f(A, A)|^{3/2}|A\cdot A|\gg |A|^{3}, ~|f(A, A)|^{3/2}|A+ A|\gg |A|^{3}.\]
Since $\min\left\lbrace |A+A|, |A\cdot A|\right\rbrace\le |A|^{\frac{9}{8}-\epsilon}$, we get $|f(A, A)|^{3/2}\gg |A|^{3-\frac{9}{8}+\epsilon}$, which concludes the proof of the theorem.
$\square$

\section{Acknowledgments}
The authors were partially supported by Swiss National Science Foundation grants 200020-162884 and 200021-175977.

\vspace{1cc}
\hfill\\
Hossein Nassajian Mojarrad\\
Department of Mathematics,\\
EPF Lausanne\\
Switzerland\\
E-mail: hossein.mojarrad@epfl.ch
\bigskip\\
Thang Pham\\
Department of Mathematics,\\
EPF Lausanne\\
Switzerland\\
E-mail: thang.pham@epfl.ch\\
\end{document}